\documentclass[preprint,12pt]{elsarticle}

\usepackage{amssymb}
\journal{Applied Mathematical Modelling}

\begin{document}

\begin{frontmatter}



\title{Numerical Modeling of Nanoparticles Transport with Two-Phase Flow in Porous Media Using Iterative Implicit Method}


%\author[M, J, S, A]{M. F. El-Amin, Jisheng Kou, Shuyu Sun, Amgad Salama}
\author[M]{M. F. El-Amin}
\author[J]{Jisheng Kou}
\author[S]{Shuyu Sun}
\author[A]{Amgad Salama}

\address[M,S,A]{Computational Transport Phenomena Laboratory (CTPL), Division of Physical Sciences and Engineering (PSE), King Abdullah University of Science and Technology (KAUST), Thuwal 23955-6900, Jeddah, Kingdom of Saudi Arabia}

\address[J]{School of Mathematics and Statistics, Hubei Engineering University, Xiaogan 432000, Hubei, China}

\begin{abstract}
In this paper, we introduce a mathematical model to describe the nanoparticles transport carried by a two-phase flow in a porous medium including gravity, capillary forces and Brownian diffusion. Nonlinear iterative IMPES scheme is used to solve the flow equation, and saturation and pressure are calculated at the current iteration step and then the transport equation is soved implicitly. Therefore, once the nanoparticles concentration is computed, the two equations of volume of the nanoparticles available on the pore surfaces and the volume of the nanoparticles entrapped in pore throats are solved implicitly. The porosity and the permeability variations are updated at each time step after each iteration loop. Two numerical examples, namely, regular heterogenous permeability and random permeability are considered. We monitor the changing of the fluid and solid properties due to adding the nanoparticles. Variation of water saturation, water pressure, nanoparticles concentration and porosity are presented graphically.   

\end{abstract}


\end{frontmatter}

% \linenumbers

\section{Introduction}
\label{Sec:Intro}
Recently, applications of nanoparticles (1-100 nm) have been reported in petroleum industry such as oil and gas exploration and production that become a promising field of research. In general, the sizes of PN are in the range of $10-500$ nm, while pore radii of a porous medium (sandstone) are from 6 to $6.3\times10^4$ nm. However, if a particle larger than a pore throat may block at the pore throat during nanoparticles transport with flow in the porous medium. There are certain types of nanoparticles, such as polysilicon nanoparticles (PN), can be used in oilfields to enhance water injection by changing the wettability through their adsorption on porous walls. The PN may be classified based on wettability of their surfaces into two types. The first type is called lipophobic and hydrophilic PN and exists in water phase only, while the second type is called hydrophobic and lipophilic PN and exists only in the oil phase. In \cite{ju1,ju2} authors have founded a mathematical model of nanoparticles transport in two-phase flow in porous media based on the formulation of fine particles transport in two-phase flow in porous media provided in Refs. \cite{liu1,liu2,liu3}. Improvements in the recovered volumes by injecting hydrophobic nanoparticles which enhance or reverse the initial reservoir wettability favoring an increase in the relative permeability of the oil phase have been reported in Ref. \cite{ony}. El-Amin et al. \cite{me1,me2} introduced a model to simulate nanoparticles transport in two-phase flow in porous media. In \cite{me3}, authors extended the model to include a negative capillary pressure and mixed relative permeabilities correlations to fit with the mixed-wet system. The model of flow and transport of nanoparticles in porous media consists mainly of five PDEs, namely, pressure, saturation, nanoparticles concentration, volume of the nanoparticles available on the pore surfaces , and finally the volume of the nanoparticles entrapped in pore throats. Moreover, the model includes variations of both porosity and permeability due to nanoparticles precipitation on the pores walls of the medium. 

The model of two-phase fluid flow in porous media is a coupled system of nonlinear time-dependent partial differential equations. Two different types of time discretization schemes are often used to solve this coupled system. The first one is the fully implicit scheme \cite{aziz,collin,daws,mont2,tan} that implicitly treats with all terms including capillary pressure. This scheme results in a system of nonlinear equations and has unconditional stability and maintains the inherent coupling of two-phase flow model. The second scheme is the IMplicit-EXplicit  (IMEX) \cite{asch,bosc,frank,hund,koto} which generally treats the linear terms implicitly and evaluates the others explicitly, and consequently. This scheme is conditionally stable, however, it has advantage that is to eliminate the nonlinearity of original equations. The IMplicit Pressure Explicit Saturation (IMPES) approach is viewed as an IMEX method, solves the pressure equation implicitly and updates the saturation explicitly. The IMPES method is conditionally stable, and hence it must take very small time step size, especially for highly heterogeneous permeable media where the capillary pressure affects substantially on the path of fluid flow. The instability of the IMPES method \cite{coats} results from the decoupling between the pressure equation and the saturation equation as well as the explicit treatment of the capillary pressure. The IMPES for two-phase flow has been improved in several versions (e.g. \cite{chen,lu,young}).

Iterative IMPES splits the equation system into a pressure and a saturation equation that are solved sequentially as IMPES \cite{lacr,blu1,blu2}. As an iterative method, the computational cost and memory required by iterative IMPES method is smaller than the fully coupled approach at each iterative step, which is more pronounced for very large size computational problems. The main disadvantage of iterative IMPES method is the decoupling of pressure and saturation equations, which results from the explicit treatment for capillary pressure. A linear approximation of capillary function is introduced to couple the implicit saturation equation into pressure equation \cite{kou2}. Kou and Sun \cite{kou1} presented an iterative version of their previous scheme proposed in \cite{kou2}. Unlike iterative IMPES, capillary pressure is not computed by the saturations at the previous iteration, but the linear approximation of capillary function at the current iteration is used, which is constructed by the saturations at the current and previous iterations. 

In this work, we use the iterative IMPES scheme introduced in \cite{kou1,kou2} to solve the flow equation of the model of nanopaticles transport in porous media. Then, we used the saturation and pressure calculated at the current iteration step to calculate the transport equation implicitly. Therefore, once we compute the nanoparticles concentration, the two equations of volume of the nanoparticles available on the pore surfaces and the volume of the nanoparticles entrapped in pore throats are solved implicitly. In our scheme, we update the porosity and the permeability variations at each time step after each iteration loop.   

\section{Modeling and Mathematical Formulation}
\subsection{Two-phase flow Model}
In this section, a mathematical model is developed to describe the nanoparticles transport carried by two-phase flow in porous media. Let us consider two-phase immiscible incompressible flow in a heterogeneous porous medium domain governed by the DarcyÕs law and the equations of mass conservation for each phase as,
\begin{equation}
\label{eq:1}
\frac{\partial \left (\phi S_{\alpha }  \right )}{\partial t} +  \nabla \cdot \mathbf{u}_{\alpha } = q_{\alpha }, \quad \alpha=w,n.
\end{equation}

\begin{equation}
\label{eq:2}
\mathbf{u}_{\alpha } = -\frac{k_{r\alpha}}{\mu_\alpha}\mathbf{K}\left ( \nabla p_\alpha + \rho_\alpha \mathbf{g} \right ), \quad \alpha=w,n.
\end{equation}

where $S_\alpha$ is the saturation, $\mathbf{u}_\alpha $ is the velocity of the phase $\alpha$. $w$ stands for the wetting phase, and $n$ stands for the nonwetting phase. $\phi$ is the porosity of the medium, $q_{\alpha }$ is the external mass flow rate. $\mathbf{K}$ is the absolute permeability tensor is chosen as $\mathbf{K}=k\mathbf{I}$, where $\mathbf{I}$ is the identity matrix and $k$ is a positive real number. $k_{r\alpha}$ is the relative permeability, $\rho_{\alpha}$ is the density, and $p_{\alpha}$ is the pressure of the phase $\alpha$, $\mathbf{g}=(0,-g)^T$ is the gravitational acceleration. $\mu_\alpha$ is the viscosity  and $k_\alpha=k_{r\alpha} \mathbf{K}$ is the effective permeability. The fluid saturations for the wetting and non-wetting are interrelated by,
\begin{equation}
\label{eq:3}
S_w + S_n = 1.
\end{equation}

Now, we describe the governing equations used in \cite{hote}, \cite{mont1} and \cite{kou1} as,
\begin{equation}
\label{eq:4}
\nabla \cdot \left(\mathbf{u}_a + \mathbf{u}_c \right) \equiv - \nabla \cdot \lambda_t \mathbf{K} \nabla {\Phi}_w -\nabla \cdot \lambda_n \mathbf{K} \nabla {\Phi}_c = q_w + q_n.
\end{equation}
and
\begin{equation}
\label{eq:5}
\frac{\partial \left (\phi S_{\alpha }  \right )}{\partial t} - q_w =  -\nabla \cdot \left(f_w \mathbf{u}_a \right) \equiv - \nabla \cdot \lambda_t \mathbf{K} \nabla {\Phi}_w. 
\end{equation}
where $f_w = \lambda_w / \lambda_t$ is the flow fraction, $\lambda_{\alpha} = k_{r\alpha} / \mu_\alpha$ is the mobility, $\Phi_w = p_w + \rho_w \mathbf{g}$ is the water pressure potential, and $\Phi_c = p_c + \left( \rho_n - \rho_w \right) \mathbf{g}$ is the capillary pressure potential. The total velocity $\mathbf{u}_t = \mathbf{u}_w + \mathbf{u}_n = \mathbf{u}_a + \mathbf{u}_c$ is defined as the sum of the two velocity variables $\mathbf{u}_a = - \lambda_t \mathbf{K} \nabla {\Phi}_w$ and $\mathbf{u}_c = - \lambda_n \mathbf{K} \nabla {\Phi}_c$. The wetting-phase velocity may be expressed by, $\mathbf{u}_w = f_w \mathbf{u}_a$. The two-phase capillary pressure can be expressed by, $p_c \left(S_w\right) = p_n - p_w$.  

\subsection{Nanoparticles Transport Model}
Ju and Fan \cite{ju1} reported that there are two types of polysillicon nanoparticles (PN) can be used in oil fields to improve oil recovery and enhance water injection, respectively. The polysillicon nanoparticles are classified based on wettability of the surface of the PN. The first type is called lipophobic and hydrophilic polysillicon nanoparticles (LHPN) and exists in water phase only, while the second type is called hydrophobic and lipophilic polysillicon nanoparticles (HLPN) and exists in the oil phase only. The sizes of PN are in the range of 10 to 500 nm, therefore, Brownian diffusion is considered. Assuming that we have a number m of size interval of the nanoparticles in water phase, the transport equation for each size interval $i$ of the nanoparticles in the water/oil phase can be written as,
\begin{equation}
\label{eq:6}
\frac{\partial \left (\phi S_\alpha C_{i,\alpha}  \right )}{\partial t} +  \mathbf{u}_\alpha \cdot \nabla C_{i,\alpha} = \nabla \cdot \left(\phi S_\alpha D_{i,\alpha} \nabla C_{i,\alpha}  \right ) - R_{i,\alpha} + Q_{i,\alpha}.
\end{equation}
where $i=1,2,É,m$. $C_{i,\alpha}$ is the volume concentrations of nanoparticles in size interval $i$ in the phase $\alpha$. 

Nanoparticles that have sizes smaller than microns have strong Brownian motion, and can bring nanoparticles very close to the pore wall. Therefore, nanoparticle may retention to decrease as flow velocity increases. The diffusion coefficient $D_i$ of the nanoparticle can be calculated using the Stokes-Einstein equation,
\begin{equation}
\label{eq:7}
D_i = \frac{k_B T}{3\pi \mu d_\rho},
\end{equation}
where $k_B$ is the Boltzmann constant, $T$ absolute temperature, $d_\rho$ is the particle diameter, and $\mu$ is the fluid viscosity. In Eq.~(\ref{eq:7}), particle diffusion constant is inversely proportional to particle diameter $d_\rho$, which increases with decreasing particle size. For instance, for water at $20^o C$ the Brownian diffusivity for a 1$-\mu m-$diameter particle is $4.3 \times 10^{-9} cm^2$ $s^{-1}$, which is small compared to solute diffusion but potentially significant over the small distances within pore spaces \cite{zhang}. $Q_{i\alpha}$ is the rate of change of particle volume belonging to a source/sink term. $R_{i\alpha}$ is the net rate of loss of nanoparticles in size interval $i$ in the phase $\alpha$. The net rate of loss of nanoparticles may be written as \cite{ju1,ju2,liu1,liu2,liu3},
\begin{equation}
\label{eq:8}
R_{i\alpha} = \frac{\partial \left( \delta \phi \right)_{i,\alpha}}{\partial t}
\end{equation}

where $\left( \delta \phi \right)_{i,\alpha} = \mathbf{v}_{1,i,\alpha} + \mathbf{v}_{2,i,\alpha}$ is the porosity change due to release or retention of nanoparticles of interval $i$ in the phase $\alpha$. $\mathbf{v}_{1,i,\alpha}$ is the volume of the nanoparticles of interval size $i$ in contact with the phase $\alpha$ available on the pore surfaces per unit bulk volume of the porous medium. $\mathbf{v}_{2,i,\alpha}$ is the volume of the nanoparticles of interval size $i$ entrapped in pore throats from the phase $\alpha$ per unit bulk volume of porous medium due to plugging and bridging. Also, $\mathbf{v}_{1,i,\alpha}$ and $\mathbf{v}_{2,i,\alpha}$ may be defined in terms of the mass of particles per unit fluid volume deposited at the pore bodies $\sigma_{1,i,\alpha}$ and pore throats $\sigma_{2,i,\alpha}$ of the porous medium as,
\begin{equation}
\label{eq:9}
\mathbf{v}_{1,i,\alpha} = \frac{\sigma_{1,i\alpha}}{\rho_b}, \quad \mathbf{v}_{2,i,\alpha} = \frac{\sigma_{2,i,\alpha}}{\rho_b}.
\end{equation}
where $\rho_b$ is the density of particulate suspensions.

At the critical velocity of the surface deposition only particle retention occurs while above it retention and entrainment of the nanoparticles take place simultaneously (Gruesbeck and Collins \cite{grues}). A modified Gruesbeck and Collins's model for the surface deposition is expressed by \cite{ju1},
\begin{equation} 
\label{eq:10}
\frac{\partial \mathbf{v}_{1,i,\alpha}}{\partial t} = \left\{\begin{array} {cc}
\gamma_{d,i,\alpha} \left \| \mathbf{u}_\alpha \right \| C_{i,\alpha},  &  \left \| \mathbf{u}_\alpha \right \| \leq  u_{\alpha c} \\ \\
\gamma_{d,i,\alpha} \left \| \mathbf{u}_\alpha \right \| C_{i,\alpha} - \gamma_{e,i,\alpha} \mathbf{v}_{1,i,\alpha} \left \| \mathbf{u}_\alpha  - u_{\alpha c}\right \|,  &  \left \| \mathbf{u}_\alpha\right \| > u_{\alpha c} 
\end{array}\right.
\end{equation}

where $\gamma_{d,i,\alpha}$ is the rate coefficients for surface retention of the nanoparticles in interval $i$ in the phase $\alpha$. $\gamma_{e,i,\alpha}$ is the rate coefficients for entrainment of the nanoparticles in interval $i$ in the phase $\alpha$. $u_{\alpha c}$ is the critical velocity for the phase $\alpha$. Similarly, the rate of entrapment of the nanoparticles in interval $i$ in the phase $\alpha$ is,
\begin{equation} 
\label{eq:11}
\frac{\partial \mathbf{v}_{2,i,\alpha}}{\partial t} = \gamma_{pt,i,\alpha} \left \| \mathbf{u}_\alpha \right \| C_{i,\alpha},  
\end{equation}
where $\gamma_{pt,i,\alpha}$ is the pore throat blocking constants.

\subsection{Porosity and Permeability Variations}
Porosity may be changed because nanoparticles deposition on the pore surfaces or blocking of pore throats. The porosity variation may be by \cite{ju1,liu1},
\begin{equation} 
\label{eq:12}
\phi = \phi_0 -   \sum_{i,\alpha} \left( \delta \phi \right)_{i,\alpha}
\end{equation}

where $\phi_0$ is the initial porosity. Also, the permeability variation due to nanoparticles deposition on the pore surfaces or blocking of pore throats may be expressed as \cite{ju2},
\begin{equation} 
\label{eq:13}
\mathbf{K} = \mathbf{K}_0  \left[ \left( 1-f \right) k_f + f \frac{\phi}{\phi_0} \right]^l
\end{equation}

where $\mathbf{K}_0$ is the initial permeability, $k_f$ is constant for fluid seepage allowed by the plugged pores. The flow efficiency factor expressing the fraction of unplugged pores available for flow is given by,
\begin{equation} 
\label{eq:14}
f = 1 -   \sum_i \gamma_{f,i} \left(\sum_{\alpha} \mathbf{v}_{2,i,\alpha} \right)
\end{equation}
$\gamma_{f,i}$ is the coefficient of flow efficiency for particles $i$. The value of the exponent $l$ has range from 2.5 to 3.5. For the nanoparticles transport carried by fluid stream in the porous media, deposition on pore surfaces and blockage in pore throats may occur. The retained particles on pore surfaces may desorb for hydrodynamic forces, and then possibly adsorb on other sites of the pore bodies or get entrapped at other pore throats.

\subsection{Initial and Boundary Conditions}
The saturation of the wetting phase in the computational domain $\Omega$ at the beginning of the flow displacing process is initially defined by,
\begin{equation} 
\label{eq:15}
S_w = S_w^0 \quad {\rm in} \quad \Omega \quad \quad {\rm at} \quad t=0.
\end{equation}
Also, the nanoparticles initial concentration of the interval $i$ in the computational domain $\Omega$ is zero, i.e.,
\begin{equation} 
\label{eq:16}
C_{i,\alpha} = 0 \quad {\rm in} \quad \Omega \quad \quad {\rm at}\quad t=0.
\end{equation}
Consequently, the initial volume of the nanoparticles of interval $i$ in contact with the phase $\alpha$ available on the pore surfaces per unit bulk volume of the porous medium is given by,
\begin{equation} 
\label{eq:17}
\mathbf{v}_{1,i,\alpha} = 0 \quad {\rm in} \quad \Omega \quad \quad {\rm at} \quad t=0.
\end{equation}
and the initial volume of the nanoparticles of interval size $i$ entrapped in pore throats from the phase $\alpha$ per unit bulk volume of porous medium due to plugging and bridging is given by,
\begin{equation} 
\label{eq:18}
\mathbf{v}_{2,i,\alpha} = 0 \quad {\rm in} \quad \Omega \quad \quad {\rm at} \quad t=0.
\end{equation}
The boundary $\partial \Omega$ of the computational domain $\Omega$ is subjected to both Dirichlet and Neumann conditions such that $\partial \Omega = \Gamma_D \cup \Gamma_N$ and $\Gamma_D \cap \Gamma_N=\o$, where $\Gamma_D$ is the Dirichlet boundary and $\Gamma_N$ is the Neumann boundary. The boundary conditions considered in this study are summarized as follow,
\begin{equation} 
\label{eq:19}
p_w \left({\rm or} \ p_n\right) = p^D \quad {\rm on} \quad \Gamma_D,
\end{equation}
\begin{equation} 
\label{eq:20}
\mathbf{u}_t \cdot \mathbf{n} = q^N \quad {\rm on} \quad \Gamma_N,
\end{equation}
where $\mathbf{n}$ is the outward unit normal vector to $\partial \Omega$, $p^D$ is the pressure on $\Gamma_D$ and $q^N$ the imposed inflow rate on $\Gamma_N$, respectively. The saturations on the boundary are subject to,
\begin{equation} 
\label{eq:21}
S_w \left({\rm or} \ S_n\right) = S^N \quad {\rm on} \quad \Gamma_N,
\end{equation}
and the concentration of the nanoparticles of interval size $i$ on the boundary is subject to,
\begin{equation} 
\label{eq:22}
C_{i,w} = C_{i,w}^0 \quad {\rm on} \quad \Gamma_N.
\end{equation}
The volume of the nanoparticles of interval size $i$ in contact with the phase $\alpha$ available on the pore surfaces per unit bulk volume of the porous medium is given by,
\begin{equation} 
\label{eq:23}
\mathbf{v}_{1,i,\alpha} = 0 \quad {\rm on} \quad \Gamma_N.
\end{equation}
and the volume of the nanoparticles of interval size $i$ entrapped in pore throats from the phase $\alpha$ per unit bulk volume of porous medium due to plugging and bridging is given by,
\begin{equation} 
\label{eq:24}
\mathbf{v}_{2,i,\alpha} = 0 \quad {\rm on} \quad \Gamma_N.
\end{equation}

\section{Iterative Implicit Method}
In this study, we consider only one interval size in the wetting phase. So, for example we drop the subscript from $C_{i,w}$, $\mathbf{v}_{1,i,\alpha}$ and $\mathbf{v}_{2,i,\alpha}$ to become $C$, $\mathbf{v}_1$ and $\mathbf{v}_2$, respectively. Similarly, we will drop the subscript $i,w,\alpha$ from all symbols including constants.
Define the time step length $\Delta t^n = t^{n+1} - t^n$, the total time interval $[0,T]$ may be divided into $N_T$ time steps as $0 = t^0< t^1  < \cdots < t^{N_T}  = T$. The current time step is represented by the superscript $n+1$. Also, the iteration loop consists of a number $N_I $ iterations of each time step. The current iteration step is denoted by $k+1$. The backward Euler time discretization is used for the equations of pressure, saturation, concentration and the two volumes to obtain,

\begin{eqnarray} 
\label{eq:25}
&& - \nabla \cdot \lambda_t \left( S_w^{n+1} \right) \mathbf{K} \left(\mathbf{v}_1^{n+1},\mathbf{v}_2^{n+1} \right) \nabla \Phi_w^{n+1} - \nabla \cdot \lambda_n \left( S_w^{n+1} \right) \ \ \ \ \ \ \ \
 \nonumber\\ 
&&
\quad \quad \quad \quad    \mathbf{K} \left(\mathbf{v}_1^{n+1},\mathbf{v}_2^{n+1} \right) \nabla \Phi_c \left( S_w^{n+1} \right) = q_w^{n+1} + q_n^{n+1},
\end{eqnarray}

\begin{eqnarray} 
\label{eq:26}
&& \phi \left(\mathbf{v}_1^{n+1}, \mathbf{v}_2^{n+1} \right) \frac{S_w^{n+1} - S_w^n}{\Delta t^n}  - q_w^{n+1} = - \nabla \cdot 
\ \ \ \ \ \ \ \
 \nonumber\\ 
&&
\ \ \ \ \lambda_t \left( S_w^{n+1} \right) \mathbf{K} \left(\mathbf{v}_1^{n+1},\mathbf{v}_2^{n+1} \right) \nabla \Phi_w^{n+1},
\end{eqnarray}

\begin{eqnarray}
\label{eq:27}
&& \phi \left(\mathbf{v}_1^{n+1},\mathbf{v}_2^{n+1} \right) \frac{S_w^{n+1} C^{n+1} - S_w^n C^n}{\Delta t^n}  + \nabla \cdot \left(\mathbf{u}_w^{n+1} C^{n+1} \right) = 
\ \ \ \ \ \ \ \
 \nonumber\\ 
&&
\quad \quad \quad \quad   
 \nabla \cdot \left(\phi \left(\mathbf{v}_1^{n+1},\mathbf{v}_2^{n+1} \right) S_w^{n+1} D \nabla C^{n+1} \right) 
\ \ \ \ \ \ \ \
 \nonumber\\ 
&&
\quad \quad \quad \quad   
  + R \left(\mathbf{u}_w^{n+1},C^{n+1},\mathbf{v}_1^{n+1} \right) + Q_w^{n+1},
\end{eqnarray}

\begin{equation}
\label{eq:28}
\frac{\mathbf{v}^{n+1}_1 - \mathbf{v}^n_1}{\Delta t^n} = \left\{\begin{array} {cc}
\gamma_d \left \| \mathbf{u}_w^{n+1} \right \| C^{n+1},  &  \left \| \mathbf{u}_w^{n+1} \right \| \leq  u_c  \\ \\
\gamma_d \left \| \mathbf{u}_w^{n+1} \right \| C^{n+1} - \gamma_e \left \| \mathbf{u}_w^{n+1} - u_c \right \| \mathbf{v}^{n+1}_1,  &  \left \| \mathbf{u}_w^{n+1} \right \| > u_c 
\end{array}\right.,
\end{equation}
and
\begin{equation}
\label{eq:29}
\frac{\mathbf{v}^{n+1}_2 - \mathbf{v}^n_2}{\Delta t^n} = \gamma_{pt} \left \| \mathbf{u}_w^{n+1} \right \| C^{n+1}.
\end{equation}
Here, because $\mathbf{u}_w^{n+1}=\mathbf{u}_w^{n+1}\left(S_w^{n+1},\Phi_w^{n+1}\right)$, we may write $R \left(\mathbf{u}_w^{n+1},C^{n+1},\mathbf{v}_1^{n+1} \right)$ as,
$$
R \left(S_w^{n+1},\Phi_w^{n+1},C^{n+1},\mathbf{v}_1^{n+1} \right) = 
$$
$$
\left\{\begin{array} {cc}
\left(\gamma_d+\gamma_{pt}\right) \left \| \mathbf{u}_w^{n+1} \right \| C^{n+1},  &  \left \| \mathbf{u}_w^{n+1} \right \| \leq  u_c  \\ \\
\left(\gamma_d + \gamma_{pt}\right) \left \| \mathbf{u}_w^{n+1} \right \| C^{n+1} - \gamma_e \left \| \mathbf{u}_w^{n+1} - u_c \right \| \mathbf{v}^{n+1}_1,  &  \left \| \mathbf{u}_w^{n+1} \right \| > u_c 
\end{array}\right..
$$

In the above equations both $\mathbf{K}$ and $\phi$ are functions in $\mathbf{v}_1$ and $\mathbf{v}_2$, that are functions of $u_w$ and $C$. On the other hand, $\mathbf{u}_w$ is a function of $\Phi_w$ and $S_w$, therefore, both $\mathbf{K}$ and $\phi$ are also functions of $\Phi_w$, $S_w$ and $C$. The system (\ref{eq:25})--(\ref{eq:29}) is fully implicit, coupled and highly nonlinear. Hence, iterative methods are often employed to solve such kind of complicated systems. In this scheme, $\mathbf{K}$ and $\phi$ will be used from the previous iteration step. Now, let us introduce the iterative formulation for the equations Eqs.~(\ref{eq:25})--(\ref{eq:29}) that is given as,
\begin{eqnarray} 
\label{eq:30}
&& - \nabla \cdot \lambda_t \left( S_w^{n+1,k} \right) \mathbf{K} \left(\mathbf{v}_1^{n+1,k},\mathbf{v}_2^{n+1,k} \right) \nabla \Phi_w^{n+1,k+1} - \nabla \cdot \lambda_n \left( S_w^{n+1,k} \right) \ \ \ \ \ \ \ \
 \nonumber\\ 
&&
\quad \quad \quad \quad\mathbf{K} \left(\mathbf{v}_1^{n+1,k},\mathbf{v}_2^{n+1,k} \right) \nabla \Phi_c \left( S_w^{n+1,k+1} \right) = q_w^{n+1} + q_n^{n+1},
\end{eqnarray}

\begin{eqnarray} 
\label{eq:31}
&& \phi \left(\mathbf{v}_1^{n+1,k}, \mathbf{v}_2^{n+1,k} \right) \frac{S_w^{n+1,k+1} - S_w^n}{\Delta t^n}  - q_w^{n+1} = - \nabla \cdot \lambda_t \left( S_w^{n+1,k} \right) \ \ \ \ \ \ \ \
 \nonumber\\ 
&&
\ \ \ \ \mathbf{K} \left(\mathbf{v}_1^{n+1,k},\mathbf{v}_2^{n+1,k} \right) \nabla \Phi_w^{n+1,k+1},
\end{eqnarray}

\begin{eqnarray}
\label{eq:32}
&& \phi \left(\mathbf{v}_1^{n+1,k},\mathbf{v}_2^{n+1,k} \right) \frac{S_w^{n+1,k+1} C^{n+1,k+1} - S_w^n C^n}{\Delta t^n}  + \nabla \cdot 
 \nonumber\\ 
&&
\left(\mathbf{u}_w^{n+1,k+1} C^{n+1,,k+1} \right) = \nabla \cdot \left\{\phi \left(\mathbf{v}_1^{n+1,k},\mathbf{v}_2^{n+1,k} \right) S_w^{n+1,k+1} \right. 
 \nonumber\\ 
&&  
\left. D \nabla C^{n+1,k+1} \right\} + R \left(S_w^{n+1,k+1},\Phi_w^{n+1,k+1},C^{n+1,k},\mathbf{v}_1^{n+1,k}\right)  
 \nonumber\\ 
&&  
+ Q_w^{n+1},
 \end{eqnarray}

\begin{eqnarray}
\label{eq:33}
&&\frac{\mathbf{v}^{n+1,k+1}_1 - \mathbf{v}^n_1}{\Delta t^n} = 
 \nonumber\\ 
&&  
\left\{\begin{array} {cc}
\gamma_d \left \| \mathbf{u}_w^{n+1,k+1} \right \| C^{n+1,k+1},  &  (c_1)  \\ \\
\gamma_d \left \| \mathbf{u}_w^{n+1,k+1} \right \| C^{n+1,k+1} - \gamma_e \left \| \mathbf{u}_w^{n+1,k+1} - u_c \right \| \mathbf{v}^{n+1,k+1}_1,  &  (c_2) 
\end{array}\right.
\end{eqnarray}
and,
\begin{equation}
\label{eq:34}
\frac{\mathbf{v}^{n+1,k+1}_2 - \mathbf{v}^n_2}{\Delta t^n} = \gamma_{pt} \left \| \mathbf{u}_w^{n+1,k+1} \right \| C^{n+1,k+1}
\end{equation}
From now upon we will refer to the condition $\left \| \mathbf{u}_w^{n+1,k+1} \right \| \leq  u_c$ by $(c_1)$ and to the condition $\left \| \mathbf{u}_w^{n+1,k+1} \right \| > u_c$ by $(c_2)$. The superscripts $k$ and $k+1$ represent the iterative steps within the current time step $n+1$. For each iteration, the variables $\lambda_w$,$\lambda_n$,$\lambda_t$ and $f_w$ is calculated using the saturation from the previous iteration. The pressure equation is solved firstly to obtain the wetting-phase pressure at the current iteration and then the DarcyÕs velocity can be calculated. Therefore, the saturation at the current iteration is computed explicitly in the current iteration. Then, the concentration and values are computed implicitly at the current time step. Finally, the permeability, porosity, and other parameters such as $\lambda_w$,$\lambda_n$,$\lambda_t$ and $f_w$ are updated. This procedure is repeated until the convergence criterion of errors has been satisfied. In this iterative scheme the capillary potential $\Phi_c$ is linearized as follows \cite{kou1},
\begin{equation} 
\label{eq:35}
\Phi_c \left(S_w^{n+1,k+1} \right) \cong \Phi_c \left(S_w^{n+1,k} \right) + {\Phi}'_c \left(S_w^{n+1,k} \right) \left[ S_w^{n+1,k+1} - S_w^{n+1,k} \right]
\end{equation}     
where ${\Phi}'_c$ is derivative of $\Phi_c$. The changes of saturation in a time step are often very small, and hence the linear approximation is reasonable. On the other hand, we use the relaxation approach to control the convergence of nonlinear iterative solvers. Therefore, the iterative scheme of pressure and saturation equations may be rewritten as,
\begin{equation} 
\label{eq:36}
- \nabla \cdot \lambda_t \left( S_w^{n+1,k} \right) \mathbf{K}^{n+1,k} \nabla \Phi_w^{n+1,k+1} - \nabla \cdot \lambda_n \left( S_w^{n+1,k} \right) \mathbf{K}^{n+1,k} \nabla \widetilde{\Phi}_c \left( \widetilde{S}_w^{n+1,k} \right) = Q^{n+1},
\end{equation}
\begin{equation} 
\label{eq:37}
\quad {\Phi}_c \left(\widetilde{S}_w^{n+1,k+1} \right) \cong \Phi_c \left(S_w^{n+1,k} \right) + {\Phi}'_c \left(S_w^{n+1,k} \right) \left[ \widetilde{S}_w^{n+1,k+1} - S_w^{n+1,k} \right],
\end{equation} 
\begin{equation} 
\label{eq:38}
\phi^{n+1,k} \frac{\widetilde{S}_w^{n+1,k+1} - S_w^n}{\Delta t^n}  + \nabla \cdot \lambda_t \left( S_w^{n+1,k} \right) \mathbf{K}^{n+1,k} \nabla \Phi_w^{n+1,k+1} = q_w^{n+1},
\end{equation}
where the relaxation equation is,
\begin{equation} 
\label{eq:39}
S_w^{n+1,k+1} = S_w^{n+1,k} + \theta_s \left( \widetilde{S}_w^{n+1,k+1} - S_w^{n+1,k} \right).
\end{equation}
where $\theta_s = \left(0,1\right]$ is a relaxation factor for saturation. In a similar manner we may write implicit iterative scheme of the equation of concentration as follow, 
\begin{eqnarray}
\label{eq:40}
%\begin{array}{c}
&& \phi^{n+1,k} \frac{S_w^{n+1,k+1} C^{n+1,k+1} - S_w^n C^n}{\Delta t^n}  + \nabla \cdot \left( \mathbf{u}_w^{n+1,k+1} C^{n+1,k+1} \right) = \nonumber\\ 
&&  
\nabla \cdot\left(\phi^{n+1,k} S_w^{n+1,k+1} D \nabla C^{n+1,k+1} \right) + 
 R \left(S_w^{n+1,k+1},\Phi_w^{n+1,k+1}, C^{n+1,k},\mathbf{v}_1^{n+1,k}\right) \nonumber\\ 
&&  + Q_c^{n+1}.
%\end{array}
\end{eqnarray}
In the above equation, we use $v_1$ from the previous iteration step in the term R. Once $C^{n+1,k+1}$ is obtained one can get $\mathbf{v}_1^{n+1,k+1}$ as follows,
\begin{equation}
\label{eq:41}
\frac{\mathbf{v}^{n+1,k+1}_1 - \mathbf{v}^n_1}{\Delta t^n} = \left\{\begin{array} {cc}
\gamma_d \left \| \mathbf{u}_w^{n+1,k+1} \right \| C^{n+1,k+1},  &  (c_1)  \\ \\
\gamma_d \left \| \mathbf{u}_w^{n+1,k+1} \right \| C^{n+1,k+1} - \gamma_e \left \| \mathbf{u}_w^{n+1,k+1} - u_c \right \| \widetilde{\mathbf{v}}^{n+1,k+1}_1,  &  (c_2) 
\end{array}\right.
\end{equation}
Finally, the equation of $\mathbf{v}_2$ may be written as,
\begin{equation}
\label{eq:42}
\frac{\mathbf{v}^{n+1,k+1}_2 - \mathbf{v}^n_2}{\Delta t^n} = \gamma_{pt} \left \| \mathbf{u}_w^{n+1,k+1} \right \| C^{n+1,k+1},
\end{equation}

\section{Numerical Tests}
In this section, we test some examples to show the performance of the presented scheme. Before presenting the numerical examples, let us introduce the necessary physical parameters used in the computations.

%\subsection{Physical parameters} 
In this study, we consider the following capillary pressure formula, 
$$
p_c = - B_c \log (S),
$$
and the normalized wetting phase saturation are correlated by, 
$$
S = \frac{S_w - S_{wr}}{1-S_{nr} - S_{wr}}, \quad 0 \leq S \leq 1,
$$
where $B_c$ is the capillary pressure parameter, $S_{wr}$ is the irreducible (minimal) water (wetting phase) saturation, and $S_{nr}$ is the residual (minimal) oil (nonwetting phase) saturation after water flooding.

Also, the expressions of the relation between the relative permeabilities and the normalized wetting phase saturation $S$ is
given as, 
$$
k_{rw} = k^0_{rw} S^a, \quad k_{rn} = k^0_{rn} \left(1 - S\right)^b,
$$
where $a$ and $b$ are positive real numbers, $k^0_{rw} = k_{rw} \left(S = 1\right)$ is the endpoint relative permeability to the wetting phase, and $k^0_{rn} = k_{rn} \left(S = 0\right)$ is the endpoint relative permeability to the non-wetting phase.

The capillary pressure function and relative permeabilities are chosen to be zero for the  residual saturations of water and oil; that is, $S=S_w$. In computation, we take the minimum of saturation as $S_{w,min}=10^{-4}$. Moreover, since the relaxation factor has some effects on the stability and efficiency of iterative methods. The choice strategy of relaxation factor \cite{kou1}, is applied in these calculations. The relaxation factor $\theta_s$ is computed where we take $\|\textbf{S}_w^{n+1,0}-\textbf{S}_w^{n+1,-1}\|=1$.  The iteration  loop continues until $\|\textbf{S}_w^{n+1,k+1}-\textbf{S}_w^{n+1,k}\|<S_{w,min}$. In our tests, we use 2-norm for vectors and matrices.

Moreover, we consider one size interval of nanoparticles suspension in the water phase. The following parameters values are used in the calculations, $\gamma_d = 16 \ m^{-1}$, $\gamma_{pt} = 1.28 \ m^{-1}$, $\gamma_e = 30 \ m^{-1}$, $\gamma_d = 16 \ m^{-1}$, $\mathbf{u}_c = 4.6 \times 10^{-6} \ m/s$, and $D = 5.6 \times10^{-8} \ m^2/s$. The nanoparticles diameter is taken as 40 nm and concentration $C_0 = 0.0$ (without nanoparticles), $0.0009, 0.004$, and $0.01$. $S_{wr}=S_{nr}=0.001$, $\phi_0=0.3$, $k_f=0.6$, $\gamma_f=0.01$. The viscosities of water and oil are 1 cP and 0.45 cP, respectively. The injection rate is  0.1 PV/year. The  relative permeabilities are quadratic, $k_{rw0}=k_{ro0}=1$, $a=b=2$, and the capillary pressure parameter is $B_c=$ 50 bar.

The domain dimension is taken as 0.3 m $\times$ 0.2 m. The computational domain is divided into 1200 uniform rectangles. The choice of relaxation factor given in Ref. \cite{kou1} is applied for the iterative method, and the three parameters are taken  as $\theta_{s,min}=0.1, \theta_{s,max}=0.9$ and $\rho=0.2$. We continue the calculation until 0.5 PVI. The time step used in this example is 0.025 day.

%\subsection{Example  1: Regular heterogenous permeabilities }

\begin{figure}\centering \includegraphics[trim=0.2cm 0.2cm 0cm 2.6cm, clip=true, width=5in,height=3in]{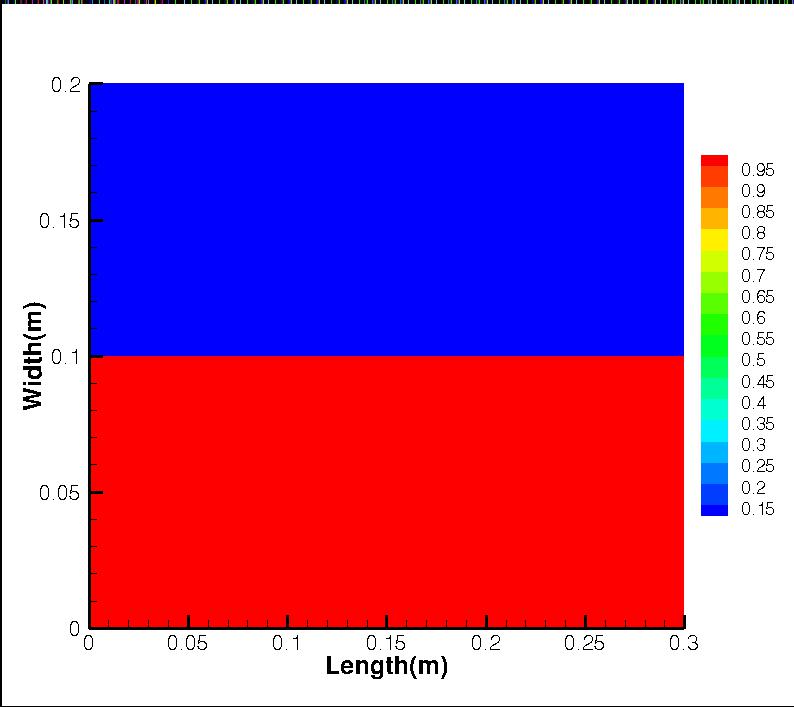} \caption{Regular heterogenous permeability (md)} \label{ex1perm} \end{figure}

\begin{figure}\centering \includegraphics[trim=0.2cm 0.2cm 0cm 2.6cm, clip=true, width=5in,height=3in]{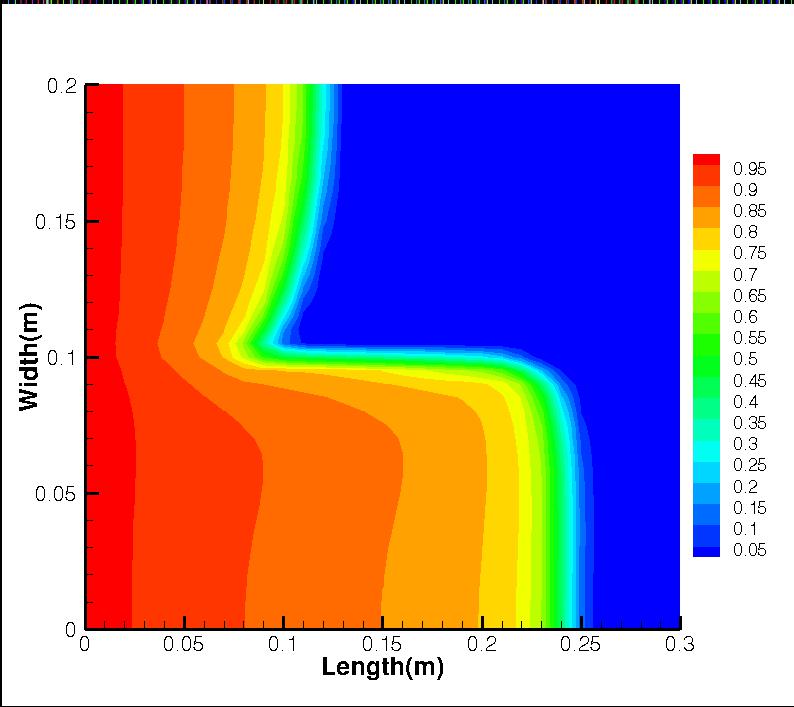} \caption{Water saturation of the regular heterogenous permeability case.} \label{ex1Sw} \end{figure}

\begin{figure}\centering \includegraphics[trim=0.2cm 0.2cm 0cm 2.6cm, clip=true, width=5in,height=3in]{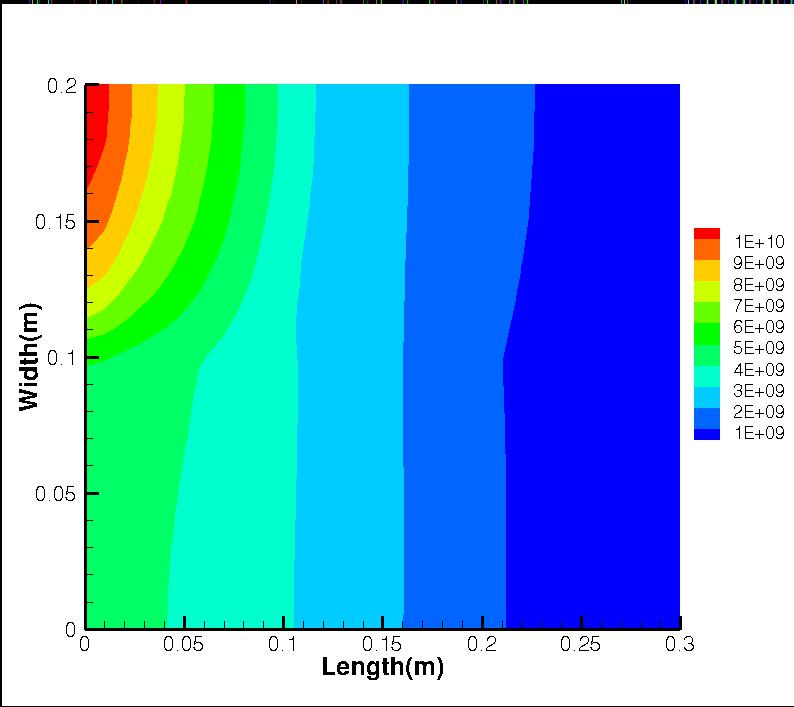} \caption{Water pressure of the regular heterogenous permeability case.} \label{ex1Pw} \end{figure}

\begin{figure}\centering \includegraphics[trim=0.2cm 0.2cm 0cm 2.6cm, clip=true, width=5in,height=3in]{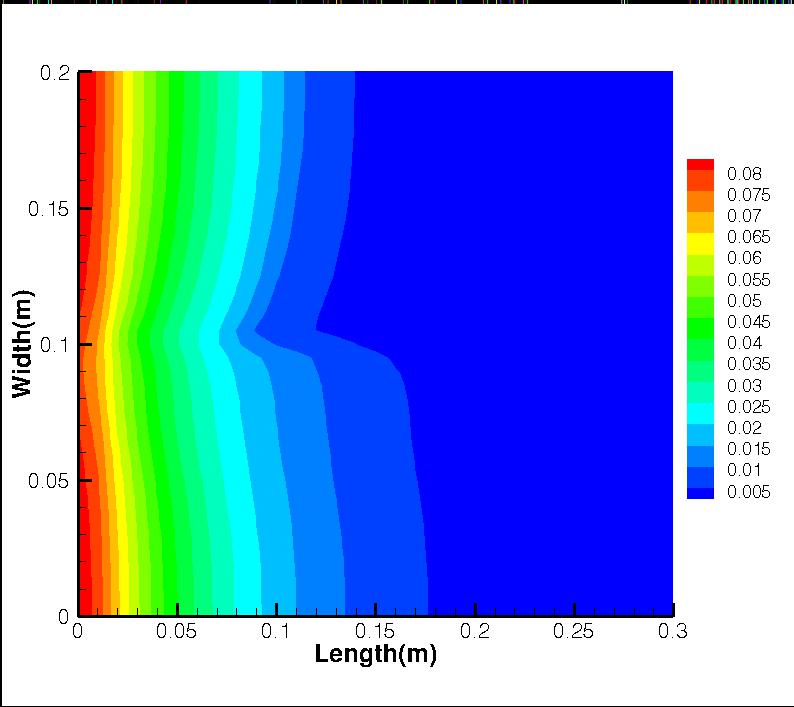} \caption{Nanoparticles concentration of the regular heterogenous permeability case.} \label{ex1C} \end{figure}

\begin{figure}\centering \includegraphics[trim=0.2cm 0.2cm 0cm 2.6cm, clip=true, width=5in,height=3in]{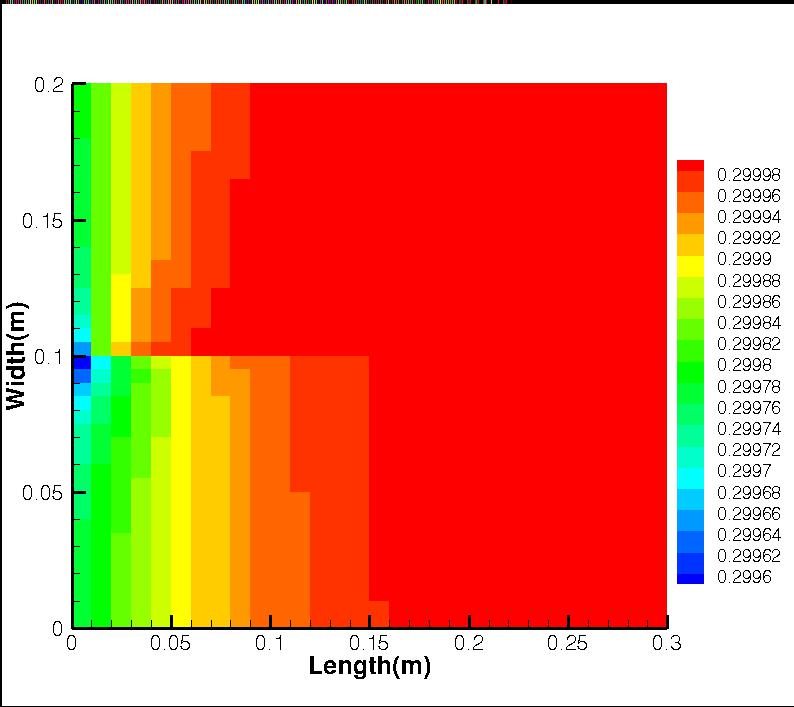} \caption{Reduced porosity of regular heterogenous permeability case.} \label{ex1poro} \end{figure}

 In the first example, the tested medium consists of  two subdomains with different configurations for the distribution of permeability as shown in Figure  \ref{ex1perm}. Figures \ref{ex1Sw},\ref{ex1Pw}, \ref{ex1C} and \ref{ex1poro} show the
distributions for water saturation, water pressure, nanoparticles concentration and porosity reduction at $0.5$ PVI, respectively.

%\begin{figure}\centering\includegraphics[width=3in,height=2in]{figures/ex1perm.eps}\caption{Heterogeneous permeabilities: Example 1. } \label{ex1perm}\end{figure}

%\subsection{Example 2:  Random permeability}

\begin{figure}\centering \includegraphics[trim=0.2cm 0.2cm 0cm 2.6cm, clip=true, width=5in,height=3in]{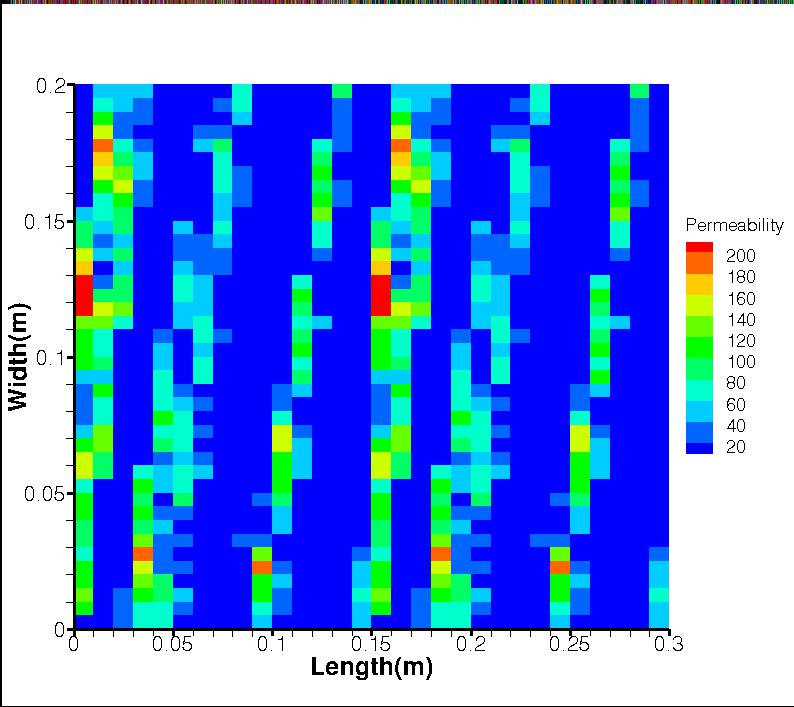} \caption{Random permeability (md)} \label{ex2perm} \end{figure}

\begin{figure}\centering \includegraphics[trim=0.2cm 0.2cm 0cm 2.6cm, clip=true, width=5in,height=3in]{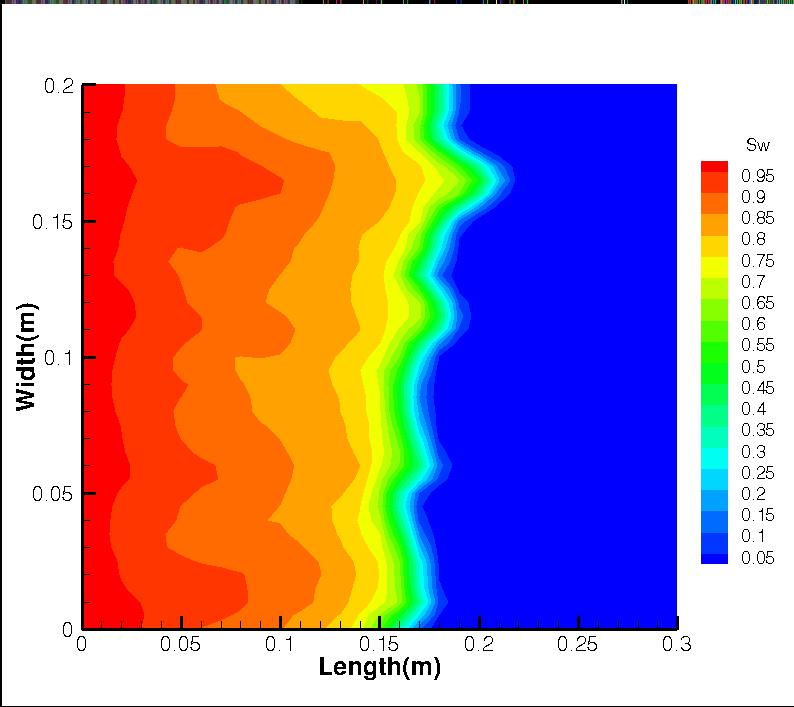} \caption{Water saturation of the random permeability case.} \label{ex2Sw} \end{figure}

\begin{figure}\centering \includegraphics[trim=0.2cm 0.2cm 0cm 2.6cm, clip=true, width=5in,height=3in]{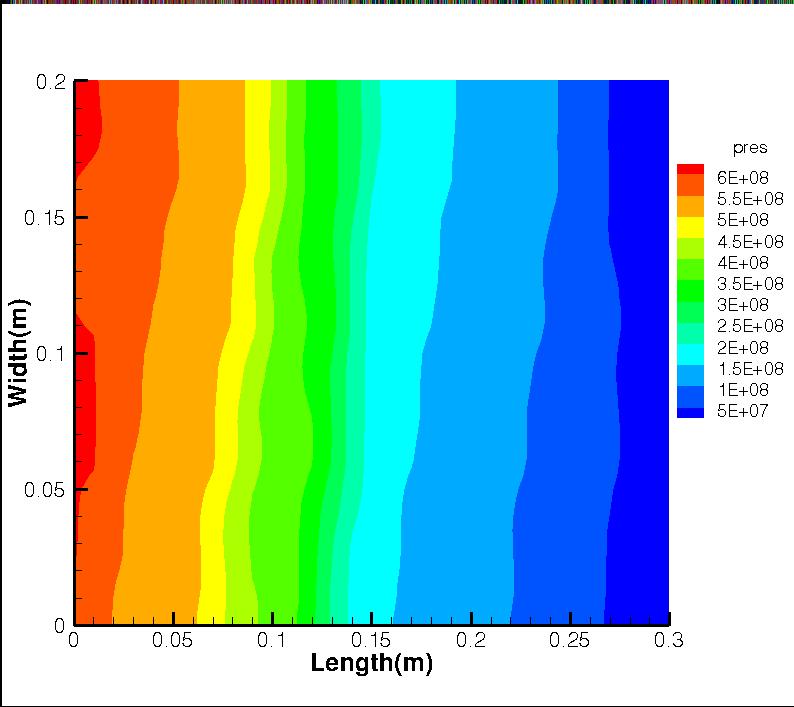} \caption{Water pressure of the random permeability case.} \label{ex2Pw} \end{figure}

\begin{figure}\centering \includegraphics[trim=0.2cm 0.2cm 0cm 2.6cm, clip=true, width=5in,height=3in]{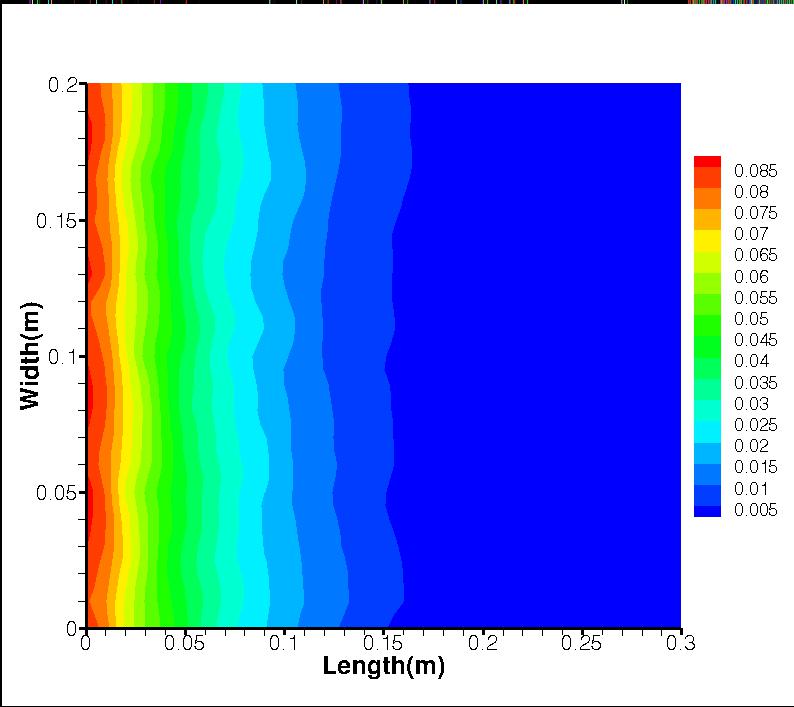} \caption{Nanoparticles concentration of the random permeability case.} \label{ex2C} \end{figure}

\begin{figure}\centering \includegraphics[trim=0.2cm 0.2cm 0cm 2.6cm, clip=true, width=5in,height=3in]{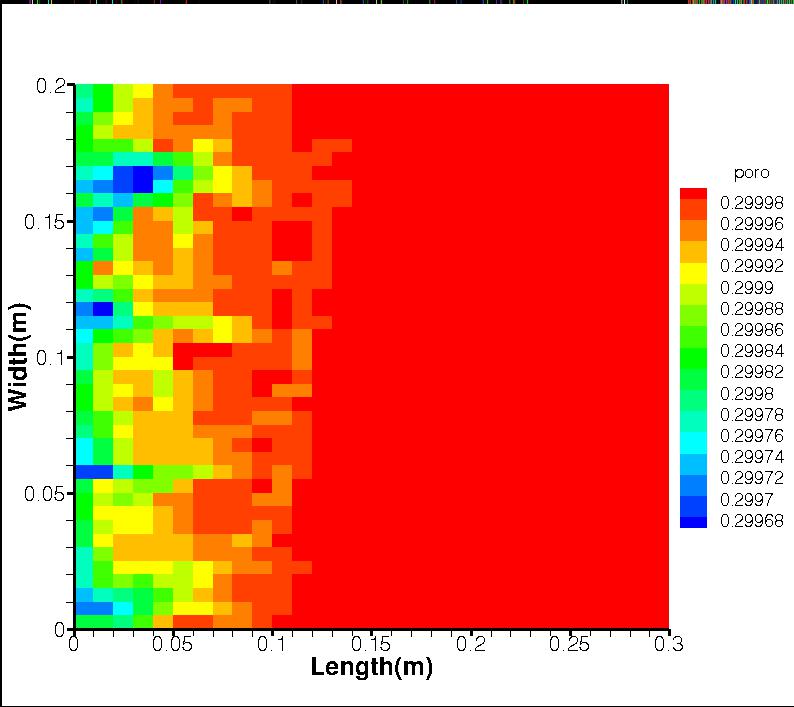} \caption{Reduced porosity of random permeability case.} \label{ex2poro} \end{figure}

 In the second example, we consider random distribution of the permeability of this medium as shown in Figure  \ref{ex2perm}. The computational results are displayed in Figures \ref{ex2Sw},\ref{ex2Pw}, \ref{ex2C} and \ref{ex2poro} show the distributions for
water saturation, water pressure, nanoparticles concentration and porosity at 0.5 PVI, respectively. We note a reduction in the porosity due to the precipitation of the nanoparticles on the porous medium walls.

\section{Acknowledgements}

This work was supported by the KAUST-UTAustin AEA project entitled: "Simulation of Subsurface Geochemical Transport and Carbon Sequestration".


\begin{thebibliography}{10}
\bibitem{hote} 
{ H.~Hoteit, and A.~Firoozabadi}, { Numerical modeling of two-phase flow in heterogeneous permeable media with different capillarity pressures}, Advances in Water Resources, 31(2008), pp.~56--73.

\bibitem{mont1}
{J.~Monteagudo and A.~Firoozabadi}, {Control-volume method for numerical simulation of two-phase immiscible flow in 2-D and 3-D discrete-fracture media}, Advances in Water Resources, 31(2004), doi: 10.1029/2003WR00299.

\bibitem{kou1} 
{J.~Kou, and S.~Sun}, { On iterative IMPES formulation for two-phase flow with capillarity in heterogeneous porous media}, Int. J. Num. Anal. Model. B, 1 (2010), pp.~20--40.

\bibitem{ju1} 
{ B.~Ju, and T.~Fan}, {Experimental study and mathematical model of nanoparticle transport in porous media}, Powder Technology, 192 (2009), pp.~195--202.

\bibitem{zhang} 
{T.~Zhang}, {Modeling of nanoparticle transport in porous media}, PhD Thesis, UT-Austin, (2012).

\bibitem{ju2} 
{B.~Ju, T.~Fan, and X.~Qiu}, {A study of wettability and permeability change caused by adsorption of nanometer structured polysilicon on the surface of porous media}, SPE--77938, SPE Asia Pacific Oil and Gas Conference and Exhibition, Melbourne, Australia, (2002).

\bibitem{liu1} 
{X.~H. Liu and F.~Civian}, {Characterization and prediction of formation damage in two-phase flow systems}, SPE--25429, Production Operations Symposium, Oklahoma City, OK, U.S.A, (1993).

\bibitem{liu2} 
{X.~H. Liu and F.~Civian}, {A multiphase mud fluid infiltration and filter cake formation model}, SPE--25215, SPE International Symposium on Oilfield Chemistry, New Orleans, LA, U.S.A., (1996).

\bibitem{liu3} 
{X.~H. Liu and F.~Civian}, {Formation damage and skin factor due to filter cake formation and fines migration in the Near- Wellbore Region}, SPE-27364, SPE Symposium on Formation Damage Control, Lafayette, Louisiana, (1994).

\bibitem{ony} 
{M.~O. Onyekonwu and N.~A. Ogolo}, {Investigating the use of nanoparticles in enhancing oil recovery} SPE-140744, Annual International Conference and Exhibition, Tinapa-Calabar, Nigeria, 2010.

\bibitem{me1} 
{M.~F. El-Amin, A.~Salama, and S.~Sun}, {Modeling and simulation of nanoparticles transport in a two-phase flow in porous media} SPE-154972, SPE International Oilfield Nanotechnology Conference and Exhibition, Noordwijk, The Netherlands, 2012.

\bibitem{me2} 
{M.~F. El-Amin, S.~Sun and A.~Salama}, {Modeling and simulation of nanoparticle transport in multiphase flows in porous media: CO$_2$ sequestration} SPE-163089, Mathematical Methods in Fluid Dynamics and Simulation of Giant Oil and Gas Reservoirs, 2012.

\bibitem{me3} 
{M.~F. El-Amin, S.~Sun and A.~Salama}, {Enhanced oil recovery by nanoparticles injection: modeling and simulation} SPE-164333, SPE Middle East Oil and Gas Show and Exhibition held in Manama, Bahrain, 10Ð13 March 2013.

\bibitem{aziz} 
{K.~Aziz and A.~Settari}, {Petroleum reservoir simulation} Applied Science Pub., London, 1979.

\bibitem{collin} 
{D.~A. Collins, L.~X. Nghiem, Y.~K. Li and J.~E. Grabenstetter}, {An efficient approach to adaptive implicit compositional simulation with an equation of state} SPE Reservoir Engineering, 7(2) (1992), pp.~259--264.

\bibitem{daws} 
{C.~N. Dawson, H.~Kl?õe, M.~F. Wheeler and C.~S. Woodward}, {A parallel, implicit, cell-centered method for two-phase flow with a preconditioned Newton- Krylov solver} Computational Geosciences, 1 (1997), pp.~215--249.

\bibitem{mont2} 
{J.~Monteagudo and A.~Firoozabadi}, {Comparison of fully implicit and IMPES formulations for simulation of water injection in fractured and unfractured media} Int. J. Numer. Meth. Engng, 69 (2007), pp.~698--728.

\bibitem{tan} 
{T.~B. Tan and N.~Kaiogerakis}, {A fully implicit, three-dimensional, three-phase simulator with automatic history-matching capability} SPE-21205,11th SPE Symposium on Reservoir Simulation, Anaheim, CA. Feb. 1991. 

\bibitem{asch} 
{U.~Ascher, S.~J. Ruuth and B.~R. Wetton}, {Implicit-Explicit methods for time-dependent partial differential equations} SIAM J. NUMER. ANAL., 32(3) (1995), pp.~797--823.

\bibitem{bosc} 
{S.~Boscarino}, {Error analysis of IMEX Runge-Kutta methods derived from differential-algebraic systems} SIAM J. Numer. Anal., 45(4) (2007), pp.~1600--1621.

\bibitem{frank} 
{J.~Frank, W.~Hundsdorfer and J.~G. Verwer}, {On the stability of implicit-explicit linear multi- step methods} Appl. Numer. Math., 25 (1997), pp.~193--205.

\bibitem{hund} 
{W.~Hundsdorfer, S.~J. Ruuth}, {IMEX extensions of linear multistep methods with general monotonicity and boundedness properties} J. Comput. Phys., 225 (2007), pp.~2016--2042.

\bibitem{koto} 
{T.~Koto}, {Stability of implicit-explicit linear multistep methods for ordinary and delay differential equations} Front. Math. China, 4(1) (2009), pp.~113--129.

\bibitem{coats} 
{K.~H. Coats}, {IMPES stability: selection of stable time steps} SPE-84924, SPE ReservoirSimulation Symposium, Houston, TX. Feb. 2001.

\bibitem{chen} 
{Z.~Chen, G.~Huan and Y.~Ma}, {Computational methods for multiphase flows in porous media} SIAM Comp. Sci. Eng., Philadelphia, 2006.

\bibitem{lu} 
{Q.~Lu}, {A parallel multiblock/multiphysics approach for multiphase flow in porous media} Ph.D Thesis, The University of Texas at Austin, 2000.

\bibitem{young} 
{L.~C. Young and R.~E. Stephenson}, {A generalized compositional approach for reservoir simulation} SPE J, 23 (1983), pp.~727--742.

\bibitem{lacr} 
{S.~Lacroix, Y.~V. Vassilevski, J.~A. Wheeler and M.~F. Wheeler}, {Iterative solution methods for modeling multiphase flow in porous media fully implicitly} SIAM J. SCI. COMPUT., 25(3) (2003), pp.~905--926.

\bibitem{blu1} 
{ B.~Lu}, {Iteratively Coupled Reservoir Simulation for Multiphase Flow in Porous Media} PhD dissertation, The University of Texas at Austin, 2008. 

\bibitem{blu2} 
{B.~Lu and M.~F. Wheeler}, {Iterative coupling reservoir simulation on high performance computers} Pet.Sci., 6 (2009), pp.43--50.

\bibitem{kou2} 
{J.~Kou and S.~Sun}, {A new treatment of capillarity to improve the stability of IMPES two-phase flow formulation} Comp. Fluids, 39 (2010), pp.~1923--1931.


\bibitem{grues} 
{C.~Gruesbeck and R.~E. Collins}, {Entrainment and deposition of fines particles in porous media}, Soc. Pet. Eng. J., 24 (1982), pp.~847--855.

\bibitem{huang} 
{D.~D. Huang, M.~M. Honarpour, and R.~Al-Hussainy}, {An improved model for relative permeability and capillary pressure incorporating wettability}, SCA, (1997), pp.~7--10.

\bibitem{skj} 
{S.~Skjaeveland, L.~Siqveland, A.~Kjosavik, W.~Hammervold, and G.~Virnovsky}, {Capillary pressure correlation for mixed-wet reservoirs}, SPE India Oil and Gas Conference and Exhibition, (1998).


\end{thebibliography}
\end{document}